\documentclass[fleqn]{mat01}
\usepackage{times,mathtimy,amssymb}
\begin{document}

\setcounter{page}{153}
\firstpage{153}

\font\xx=msam5 at 10pt
\def\qed{\mbox{\xx{\char'03\!}}}

\makeatletter
\def\artpath#1{\def\@artpath{#1}}
\makeatother
\artpath{c:/pm2140}

\newcommand{\psp}{\mathbf{P}}
\newcommand{\Fq}{\mathbf{F}_q}
\newcommand{\Ov}{\mathbf{O}_v}
\newcommand{\Oh}{\mathbf{O}}
\newcommand{\A}{\mathbf{A}}
\newcommand{\Gv}{\Gamma_v}
\newcommand{\bu}{\bullet}
\newcommand{\fbu}{F_\bu}
\newcommand{\h}{\mathbf{H}}
\newcommand{\kappabar}{\overline{\kappa}}
\newcommand{\bsl}{\backslash}
\newtheorem{theor}{\bf Theorem}
\newtheorem{cor}{\rm COROLLARY}
\newtheorem{propo}{\rm PROPOSITION}
\newtheorem{lem}{Lemma}

\newcommand{\R}{\mathbb{R}}

\title{Reduction theory for a rational function field}

\markboth{Amritanshu Prasad}{Reduction theory for a rational function field}

\author{AMRITANSHU PRASAD}

\address{Max-Planck-Institut f\"{u}r Mathematik,
 Postfach 7280, D-53072 Bonn, Germany}

\volume{113}

\mon{May}

\parts{2}

\Date{MS received 9 September 2002; revised 16 October 2002}

\begin{abstract}
Let $G$ be a split reductive group over a finite field $\Fq$. Let
$F=\Fq(t)$ and let $\A$ denote the ad\`eles of $F$. We show that every
double coset in $G(F)\bsl G(\A)/ K$ has a representative in a maximal
split torus of $G$. Here $K$ is the set of integral ad\`elic points of
$G$. When $G$ ranges over general linear groups this is equivalent to
the assertion that any algebraic vector bundle over the projective line
is isomorphic to a direct sum of line bundles.
\end{abstract}


\maketitle

\section{Introduction}

Let $F$ be a global field, $\A$ its ring of ad\`eles and $G$ a reductive
group defined over $F$. The theory of automorphic forms involves the
study of spaces of functions on $G(F)\bsl G(\A)$ as representations of
$G(\A)$. The functions involved are often required to be right invariant
under certain large compact subgroups $K$ of $G(\A)$ because (among
other reasons) the double coset space $G(F)\bsl G(\A)/K$ admits nice
interpretations. For example, the classical study of the upper half
plane modulo the action of arithmetic subgroups of the real special
linear group is a special case of the above when $F$ is the field of
rational numbers (see e.g., (\cite{PS}, {\S1}). Another special case, which
corresponds to taking $F$ to be a field of rational functions in one
variable and $G$ to be $GL(2)$ is discussed by Weil in \cite{W}. When
$F$ is a function field, Harder describes a fundamental domain for the
action of $G(F)$ on $G(\A)$ in (\cite{H74}, {\S1}) using results from
\cite{H68} and \cite{H69}. This is an analogue of the Seigel domain
described by Godement in \cite{God} for $F=\mathbf{Q}$. Proposition
\ref{prop:Omega} in this article is analogous to these results and the
proof proceeds along the lines of \cite{God}. Harder's description of
the fundamental domain is a very basic result in the theory of
automorphic forms over function fields (see e.g., \cite{Laumon2}, \S9 and
Appendix E).

From now on let $G$ be a split reductive group defined over a finite
field $\Fq$ with $q$ elements. Fix a Borel subgroup $B$ defined over
$\Fq$ with unipotent radical $N$, and a maximal $\Fq$-split torus $T$
contained in $B$. Set $F=\Fq(t)$. For a valuation $v$ of $F$, we denote
the corresponding local field by $F_v$ and its ring of integers by
$\Ov$. For each $v$, fix a uniformizing element $\pi_v\in F\cap \Ov$. In
particular, fix $\pi_\infty=t^{-1}$ as a uniformizing element at the
place $\infty$ whose local field is $\Fq((t^{-1}))$. Let $K$ be the
maximal compact subgroup $\prod_v G(\Ov)$ of $G(\A)$. This article
concerns the double coset space
\begin{equation*}
G(F)\bsl G(\A) / K
\end{equation*}
which may be interpreted as the set of isomorphism classes of principal
$G$-bundles on the projective line. In \cite{Groth}, Grothendieck proves
that when $G$ is a complex reductive group any holomorphic $G$-bundle
over the complex projective line admits a reduction of structure group
to a maximal torus. (In fact this result has been attributed to Dedekind
and Weber for $G=GL(n)$ by Geyer (\cite{Geyer}, \S6) who deduces it from a
statement in (\cite{DW}, \S22).) In our ad\`elic setting, this should
correspond to the assertion that every double coset has a representative
in $T(\A)$.

Let $X_*(T)$ denote the lattice $\mathrm{Hom} (\mathbf{G}_m,T)$ of
algebraic co-characters of $T$. Given $\eta\in X_*(T)$, and a valuation
$v$ denote by $\pi_v^\eta$ the element $\eta(\pi_v)\in T(F_v)\subset
T(\A)$. Recall that $\eta\in X_*(T)$ is called \emph{antidominant} if
$|\alpha_i\circ\eta(\pi_v)|_v\geq 1$ for each simple root $\alpha_i$
(see \S\ref{section:fundamental_representations}). Precisely stated, the
main result of this article is the following:

\begin{theor}[\!]
\label{theorem:global_reduction}
Every double coset in
\begin{equation*}
G(F)\backslash G(\A) / K
\end{equation*}
has a unique representative of the form $(t^{-1})^\eta${\rm ,} where $\eta \in
X_*(T)$ is antidominant.
\end{theor}

In \S\ref{sec:global_proof}, we will deduce Theorem
\ref{theorem:global_reduction} from the following local result which is
proved in \S \ref{sec:local_proof}. Let $\fbu$ be the local field
$\Fq((\pi))$ of Laurent series in $\pi$ with coefficients in $\Fq$. It
contains, as its ring of integers, the discrete valuation ring
$\Oh=\Fq[[\pi]]$, and as a discrete subring, the polynomial ring
$R=\Fq[\pi^{-1}]$. Let $\Gamma=G(R)$.

\begin{theor}[\!]
\label{theorem:local}
Every double coset in
\begin{equation*}
\Gamma \backslash G(F_\bu) / G(\Oh)
\end{equation*}
has a unique representative of the form $\pi^{\eta},$ where $\eta \in
X_*(T)$ is antidominant.
\end{theor}

The main results proved in this article should be known to the experts,
but we have not found them in the literature beyond the case of $GL(2)$,
for which Theorem \ref{theorem:local} is proved in (\cite{W}, \S3). The
results proved in this paper have played an important role in the
author's work \cite{P02}, as well as
in the work of other authors on $\Fq(t)$ \cite{E,Ans,KR}.

\section{Normed local vector spaces}

Let $V$ be a vector space defined over $\Fq$. Let $\mathbf{e}_1,
\ldots,\mathbf{e}_n$ be a basis of the free $\Oh$-module $V(\Oh)$ (so
that $V(\Oh)$ is isomorphic to the free $\Oh$-module generated by the
$\mathbf{e}_i$s). Given a vector $\mathbf{x}\in V(F_\bu)$, we may write
$\mathbf{x}=x_1 \mathbf{e}_1+\cdots+x_n\mathbf{e}_n$, uniquely, with
$x_i\in F_\bu$.\break Define
\begin{equation}
\label{eq:defn_of_norm}
\|\mathbf{x}\|=\sup\{|x_1|,\ldots , |x_n|\}.
\end{equation}
\setcounter{lem}{3}
\begin{lem}\label{lemma:canonicity}
If $g\in GL(V(\Oh)),$ then $\|\mathbf{x} g\| = \| \mathbf{x} \|$.
\end{lem}
\begin{proof}
Let $(g_{ij})$ be the matrix of $G$ with respect to the basis chosen
above. Let $\mathbf{y}=\mathbf{x}g$. If $\mathbf{y}=y_1
\mathbf{e}_1+\cdots + y_n \mathbf{e}_n$, then
\begin{equation*}
y_j = \sum_{i=1}^n x_i g_{ij}
\end{equation*}\pagebreak

\noindent and\vspace{-1pc}
\begin{equation*}
\begin{array}{rclr}
\| \mathbf{y} \| & = & \sup\limits_{1\leq j \leq n} | \sum\limits_{i=1}^n x_i g_{ij} | & \\[1pc]
& \leq & \sup\limits_{1\leq j \leq n}\  \sup\limits_{1\leq i \leq n} |x_i g_{ij}| & \mbox{ (ultrametric inequality)}\\[1pc]
& \leq & \sup\limits_{1\leq j \leq n}\  \sup\limits_{1\leq i \leq n} |x_i| & \mbox{ (since } g_{ij} \in \Oh \mbox{)}\\[1pc]
& = & \| \mathbf{x} \|.
\end{array}
\end{equation*}
Hence
\begin{equation*}
\| \mathbf{y}\| \leq \| \mathbf{x}\|.
\end{equation*}
We may apply the same reasoning to $g^{-1}$ to show that
\begin{equation*}
\| \mathbf{x}\| \leq \| \mathbf{y}\|.
\end{equation*}
Therefore,
\begin{equation*}
\| \mathbf{y}\| = \| \mathbf{x}\|.
\end{equation*}
\hfill \qed
\end{proof}

\setcounter{cor}{4}
\begin{cor}
\label{cor:canonicity}
$\left.\right.$\vspace{.5pc}

\noindent The norm $\|\cdot \|$ is independent of our choice of basis of $V(\Oh)$.
\end{cor}

\begin{proof}
The coordinates of a vector with respect to two different bases differ
by a matrix with entries in $\Oh$. The argument in the proof of Lemma
\ref{lemma:canonicity} shows that the norms with respect to two
different bases are equal.\hfill \qed
\end{proof}

\setcounter{lem}{5}
\begin{lem}\label{lemma:ultrametric}
The norm $\|\cdot \|$ satisfies the {\rm ultrametric triangle
inequality,} i.e.{\rm ,}  for vectors $\mathbf{x}$, $\mathbf{y}$ in
$V(F_\bu),$
\begin{equation*}
\| \mathbf{x} + \mathbf{y} \| \leq \sup\{ \| \mathbf{x}\|,\|\mathbf{y}\|
\}.
\end{equation*}
\end{lem}

\begin{proof}
Write $\mathbf{x}=x_1\mathbf{e}_1+\cdots + x_n \mathbf{e}_n$ and
$\mathbf{y} = y_1\mathbf{e}_1+\cdots + y_n \mathbf{e}_n$. Then
\begin{align*}
\| \mathbf{x}+\mathbf{y} \| & = \sup\{|x_1+y_1|, \ldots, |x_n+y_n|\}\\[.3pc]
& \leq \sup\{\sup\{|x_1|,|y_1|\}, \ldots, \sup\{|x_n|,|y_n|\}\}\\[.3pc]
& = \sup\{|x_1|,|y_1|, \ldots, |x_n|,|y_n|\}\\[.3pc]
& = \sup\{ \|\mathbf{x}\|,\|\mathbf{y}\| \}.
\end{align*}
\hfill \qed
\end{proof}

\begin{lem}
For a scalar $\lambda \in F_\bu$ and a vector $\mathbf{x}\in V(F_\bu),$
\begin{equation*}
\| \lambda \mathbf{x}\| = |\lambda| \|\mathbf{x}\|.
\end{equation*}
\end{lem}
\begin{lem}\label{lemma:bound}
If $g\in GL(V(F_\bu)),$ then there is a constant $C_g>0,$ such that for
any vector $\mathbf{x}\in V(F_\bu),$
\begin{equation*}
\|\mathbf{x}g\| \leq C_g \| \mathbf{x}\|.
\end{equation*}
\end{lem}

\begin{proof}
Suppose that $g$ has matrix $(g_{ij})$, and $\mathbf{x}$ has coordinates
$(x_1, \ldots, x_n)$ with respect to the basis $\mathbf{e}_1,\ldots,
\mathbf{e}_n$. Then
\begin{align*}
\|\mathbf{x}g\| & = \sup\left\{\left|\sum_{i=1}^nx_i g_{i1}\right|,\ldots,
\left|\sum_{i=1}^nx_i g_{in}\right|\right\}\\
& \leq \sup_{1\leq j \leq n} \sup_{1\leq i \leq n} |g_{ij}|
\|\mathbf{x}\|.
\end{align*}
Therefore, let
\begin{equation*}
C_g=\sup_{1\leq j \leq n} \sup_{1\leq i \leq n} |g_{ij}|.
\end{equation*}
\hfill \qed
\end{proof}

\begin{lem}\label{lemma:boundone}
If $\mathbf{x}\in V(R)$ is a non-zero vector then $\|\mathbf{x}\|\geq 1$.
\end{lem}

\begin{proof}
By Corollary \ref{cor:canonicity}, we may assume that the elements
$\mathbf{e}_i$ of a basis used to define $\|\cdot \|$ lie in $V(\Fq)$.
Then at least one coordinate of $\mathbf{x}$ is non-zero in $R$. But any
non-zero element in $R$ has norm at least one. Therefore,
$\|\mathbf{x}\|\geq 1$.\hfill \qed
\end{proof}
\setcounter{propo}{9}
\begin{propo}\label{prop:lower_bound}$\left.\right.$\vspace{.5pc}

\noindent For any non-zero vector $\mathbf{x} \in V(\Fq)$ and any $g\in
GL(V(F_\bu)),$ there is a positive constant $E$ such that for all
$\gamma \in GL(V(R)),$
\begin{equation*}
\| \mathbf{x} \gamma g \| \geq E.
\end{equation*}
Consequently{\rm ,} for any subset $S$ of $GL(V(R)),$ the set $\{\|\mathbf{x}
sg\|: s\in S\}$ has a positive minimal element.
\end{propo}

\begin{proof}
Applying Lemma \ref{lemma:bound} to $g^{-1}$, and Lemma
\ref{lemma:boundone} to $\mathbf{x}\gamma$ (which lies in $V(R)$), we
have
\begin{equation*}
\|\mathbf{x}\gamma g \| \geq C_{g^{-1}}\| \mathbf{x}\gamma\| \geq C_{g^{-
1}}>0.
\end{equation*}
The second part of the assertion follows by noting that the values taken
by the norm $\|\cdot\|$ are of the form $q^j$, where $j$ is an
integer.\hfill \qed
\end{proof}

\section{Fundamental representations}
\label{section:fundamental_representations}

Let $\alpha_1,\ldots, \alpha_r$ be the simple roots with respect to $B$
in the root system $\Phi(G,T)$ of $G$ with respect to $T$. Let
$W=N_G(T)/T$ be the Weyl group of $G$ with respect to $T$. To each
simple root $\alpha_i$, we associate an element $s_i$ of order two in
$W$ in the usual way.

Given a subset $D$ of $\{1,\dots,r\}$, let $W_D$ denote the subgroup of
$W$ generated by $\{s_j| j\in D\}$, and let $P_D$ denote the parabolic
subgroup $BW_DB$ of $G$ containing $B$. This group has a Levi
decomposition
\begin{equation*}
P_D=L_DU_D,
\end{equation*}
where $L_D$ is a reductive group of rank $|D|$ and $U_D$ is the
unipotent radical of $P_D$. $L_D\cap B$ is a Borel subgroup for $L_D$
containing the split torus $T$. The set of simple roots of $L_D$ with
respect to $L_D\cap B$ is $\{ \alpha_j | j \in D \}$. Denote by $P_i$
(resp., $L_i$, $U_i$) the parabolic subgroup (resp., Levi subgroup,
unipotent subgroup) corresponding to the set
$\{1,\ldots,i-1,i+1,\ldots,r\}$. These are the maximal proper parabolic
subgroups of $G$ containing $B$.
\setcounter{theor}{10}
\begin{theor}[\!]{\bf \cite{Chevalley}}
\label{theorem:godement}
There exist irreducible finite dimensional representations
$(\rho_i,V_i)$ of $G,$ vectors $\mathbf{v}_i\in V_i(\Fq)$ that are
unique up to scaling{\rm ,} and characters $\Delta_i:P_i \to \mathbf{G}_m,$
for $i=1,\ldots,r$ all defined over $\Fq,$ such that\vspace{-.5pc}
\end{theor}
\begin{enumerate}
\renewcommand{\labelenumi}{\arabic{enumi}.}
\item \label{item:delta}{\it $P_i$ is the stabilizer of the line generated by
$\mathbf{v}_i$ and $\mathbf{v}_i \rho_i(p) = \Delta_i(p) \mathbf{v}_i$
for each $p\in P_i$ for $i=1,\ldots,r$.}
\item \label{item:mu}{\it The restrictions $\mu_i$ to $T$ of $\Delta_i$s are
antidominant weights of $T$ with respect to $B,$ which generate
$X^*(T)\otimes \mathbf{Q}$ as a vector space over the rational numbers.}
\end{enumerate}

\section{Ordering by roots}
\label{sec:ordering}
\setcounter{lem}{11}
\begin{lem}\label{lemma:lifting_characters}
Let $L$ be a Levi subgroup of $G$ associated to a parabolic subgroup $P$
containing $B$. Then there is a canonical surjection
\begin{equation*}
G(\fbu)/G(\Oh) \stackrel{^{\Phi^G_L}}{\longrightarrow} L(\fbu)/L(\Oh).
\end{equation*}
\looseness 1 If $Q=MN$ is a parabolic subgroup of $G$ containing $B$ and contained in
$P,$ then $M$ is a Levi subgroup for $L$ corresponding to the parabolic
subgroup $L\cap Q$ of $L,$ and $\Phi^L_M\circ\Phi^G_L=\Phi^G_M$.
\end{lem}

\begin{proof}
Given $g\in G(\fbu)$, we may use the Iwasawa decomposition to write
$g=luk$, where $l\in L(\fbu)$, $u\in U(\fbu)$ and $k\in G(\Oh)$.
Moreover, if $g=l'u'k'$ is another such decomposition, then, setting
$l_0=l^{\prime-1}l$ and $k_0=k'k^{-1}$,
\begin{equation*}
u'{^{-1}}l_0u=k_0 \in G(\Oh).
\end{equation*}
On the other hand,
\begin{equation*}
k_0 = u'{^{-1}}l_0u = l_0l_0^{-1}u'{^{-1}}l_0u.
\end{equation*}
Since $L$ normalizes $U$, $l_0^{-1}u'{^{-1}}l_0 \in U(\fbu)$, and hence,
setting $u_0=l_0^{-1}u'{^{-1}}l_0u\in\break U(\fbu)$,
\begin{equation*}
l_0=k_0u_0\in G(\Oh) U(F_\bu) \cap L(F_\bu).
\end{equation*}
Therefore $l_0u_0^{-1}=k_0\in G(\Oh)\cap P(F_\bu)=P(\Oh)$, so that
$l_0\in L(\Oh)$. This shows that $luk\mapsto l$ induces a well defined
map $\Phi^G_L: G(\fbu)/G(\Oh) \to L(\fbu)/L(\Oh)$. It is clear that this
map is surjective. To see that $\Phi^L_M\circ\Phi^G_L=\Phi^G_M$, note
that we may write $g=muk$ with $m\in M(\fbu)$, $u\in N(\fbu)$ and $k\in
G(\Oh)$. But $N(\fbu)=(N(\fbu)\cap L(\fbu))U(\fbu)$, so we may write
$u=u_1u_2$, where $u_1\in N(\fbu)\cap L(\fbu)$ and $u_2\in U(\fbu)$.
Therefore, we see that $mM(\Oh)=\Phi^L_M(mu_1)=\Phi^G_M(g)$.\hfill \qed
\end{proof}

In the sequel we denote $\Phi^G_T$ simply by $\Phi$. Define
\begin{equation}
\label{eq:defn_Omega}
\Omega_G := \{ g\in G(\fbu)\; : \;|\alpha_i\circ\Phi(g)| \geq 1 \mbox{ for } i=1,\ldots,r\}.
\end{equation}

\setcounter{propo}{13}
\begin{propo}
\label{prop:Omega}$\left.\right.$\vspace{.5pc}

\noindent $G(\fbu)=\Gamma \Omega_G$.
\end{propo}

\begin{proof}$\left.\right.$\vspace{.5pc}

\noindent{\it The rank one case} (following \cite{W}): Here $G$ has one simple
root $\alpha_1$, and one fundamental representation $(\rho_1,V_1)$ and a
vector $\mathbf{v}_1\in V_1(\Fq)$ such that for any element $p$ in the
parabolic subgroup $B=TN$, where $N$ is the unipotent radical of $B$,
\begin{equation}
\label{eq:mu1}
\mathbf{v}_1\rho_1(b)=\Delta_1(b) \mathbf{v}_1,
\end{equation}
where the character $\Delta_1:B\mapsto \mathbf{G}_m$ (defined over
$\Fq$) restricts to an anti-dominant weight $\mu_1$ on the maximal split
torus $T$. Let $g\in G(F_\bu)$. We wish to show that $g\in
\Gamma\Omega_G$. To this end, by Proposition \ref{prop:lower_bound}, and by
replacing $g$, if necessary by an appropriate element of $\Gamma g$, we
may assume that $g$ has the property that
\begin{equation}
\label{eq:assumption}
\|\mathbf{v}_1 \rho_1(\gamma g)\|\geq \|\mathbf{v}_1 \rho_1(g)\|\quad \mbox{
for all } \gamma \in \Gamma.
\end{equation}
Write $g = tnk$, where $t\in T(F_\bu)$, $n\in N(F_\bu)$ and $k\in
G(\Oh)$. By Theorem \ref{theorem:godement} and Lemma~\ref{lemma:canonicity},
\begin{equation}
\label{eq:resumption}
\| \mathbf{v}_1 \rho(g) \| = |\Delta_1(t)|\|\mathbf{v}_1\|=|\mu_1(t)|.
\end{equation}
Fix an isomorphism $u_{\alpha_1}:\mathbf{G}_a\to N$ defined over $\Fq$,
and let $x\in \fbu$ be such that $n=u_{\alpha_1}(x)$.
Choose $\sigma$ in the nontrivial $T(\Fq)$-coset of $N_GT(\Fq)$.
Note that if $S\in R$, then $\sigma u_{\alpha_1}(S)\in \Gamma$,
therefore, using Proposition \ref{prop:lower_bound},
\begin{align*}
|\mu_1(t)|& = \|\mathbf{v}_1\rho_1(g)\|\\
& \leq \| \mathbf{v}_1 \rho_1(\sigma u_{\alpha_1}(S)t
u_{\alpha_1}(x))\|\\
& = \| \mathbf{v}_1 \rho_1(^\sigma t \sigma u_{\alpha_1}(\alpha_1(t)^{-
1}(S+\alpha_1(t)x)))\|\\
& = |\mu_1(t)|^{-1}\|\mathbf{v}_1\rho_1(u_{-\alpha_1}(\alpha(t)^{-
1}S+x))\|.
\end{align*}
Here $u_{-\alpha_1}=\sigma u_{\alpha_1} \sigma^{-1}$, and its image is
the root subgroup for $-\alpha_1$. The element
$u_{-\alpha_1}(\alpha(t)^{-1}S+x)$ lies in the derived group of $G$
which is isomorphic to either $SL_2$ or $PGL_2$ in the rank one case.
When the derived group of $G$ is isomorphic to $SL_2$, we may take $V_1$
to be the right action of $SL_2$ on the space of $1\times 2$-matrices by
right multiplication. One may take the torus $T$ to consist of diagonal
matrices in $SL_2$, $B$ the upper triangular matrices in $SL_2$ and
$\mathbf{v}_1$ to be the vector $(0,1)$. Calculating with matrices, one
may verify that
\begin{equation*}
\|\mathbf{v}_1\rho_1(u_{-\alpha_1}(\alpha(t)^{-1}S+x))\|\leq \sup\{
1,|\alpha(t)^{-1}S+x| \}.
\end{equation*}
Therefore,
\begin{equation}
\label{eq:ineq}
\sup\{ 1, |\alpha_1(t)^{-1}S+x|\}\geq |\mu_1(t)|^2.
\end{equation}
Choose $S$ in $R$ such that $|S+\alpha(t)x|<1$. Then
$|\alpha_1(t)^{-1}S+x|<|\alpha_1(t)|^{-1}$. Suppose that
$|\alpha_1(t)^{-1}S+x|\geq|\mu_1(t)|^2$. Then
$|\alpha_1(t)|^{-1}>|\mu_1(t)|^2$. This is impossible, since
$\alpha_1(t)^{-1}=\mu_1(t)^2$. It follows that
$|\alpha_1(t)^{-1}S+x|<|\mu_1(t)|^2$. Therefore, (\ref{eq:ineq}) can
hold only if $1\geq |\mu_1(t)|^2$, which is the same as
$|\alpha_1(t)|\geq 1$. This completes the proof of Proposition
\ref{prop:Omega} when the derived group of $G$ is isomorphic to $SL_2$.

When the derived group of $G$ is isomorphic to $PGL_2$, then $G$ is the
product of its centre with $PGL_2$. Therefore, the assertion of
Proposition \ref{prop:Omega} for $G$ follows from that for $PGL_2$.
However, the assertion for $PGL_2$ follows easily from that for $GL_2$.
The derived group of $GL_2$ is $SL_2$, hence the proposition holds for
$GL_2$ by the argument in the previous paragraph, completing the proof
of Proposition \ref{prop:Omega} in the rank one case.\vspace{.5pc}

\noindent {\it The general case}:
Let $G$ be a group of rank $r$, and $g\in G(F_\bu)$.
By modifying $g$ on the left by an element of $\Gamma$, we may, for the
purposes of this proof, assume, using the second assertion of
Proposition \ref{prop:lower_bound}, that
\begin{equation}
\label{eq:first_inequality}
\|\mathbf{v}_1\rho_1(g)\|\leq \|\mathbf{v}_1 \rho_1(\gamma g)\|\quad \mbox{
for all } \gamma \in \Gamma.
\end{equation}
Note that if $\gamma\in P_1(\fbu)\cap\Gamma$, then
$\mathbf{v}_1\rho_1(\gamma g)=\Delta_1(\gamma)\mathbf{v}_1\rho_1(g)$.
Since $\Delta_1(\gamma)\in \Fq[\pi^{-1}]^\times$,
$|\Delta_1(\gamma)|=1$. Therefore, $\| \mathbf{v}_1\rho_1(\gamma g)\|
=\|\Delta_1(\gamma)\mathbf{v}_1\rho_1(g)\|$. We may use the second
assertion of Proposition \ref{prop:lower_bound} again, to assume, for
the purposes of this proof, that
\begin{equation}
\|\mathbf{v}_2\rho_2(g)\|\leq \|\mathbf{v}_2 \rho_2(\gamma g)\|\quad \mbox{
for all } \gamma \in \Gamma\cap P_1(\fbu)
\end{equation}
while preserving (\ref{eq:first_inequality}).
Continuing in this manner, we may assume that
\begin{equation}
\label{eq:jth_inequality}
\|\mathbf{v}_j\rho_j(g)\|\leq \|\mathbf{v}_j\rho_j(\gamma g)\|\quad \mbox{
for all } \gamma \in \Gamma\cap P_1(F)\cap \ldots \cap P_{j-
1}(F),
\end{equation}
for $j=1,\ldots,r$. Therefore, it suffices to prove the following:
\setcounter{lem}{21}
\begin{lem}
\label{lemma:otherthing}
If an element $g\in G(\fbu)$ satisfies the
inequalities $(\ref{eq:jth_inequality})$ for each integer $1\leq j \leq r,$
then $g\in \Omega_G$.
\end{lem}

The proof of Proposition
\ref{prop:Omega} in the rank one case shows that Lemma
\ref{lemma:otherthing} is true when $G$ is of semisimple rank one. We
prove it in general assuming the validity of Theorem
\ref{theorem:local} in the rank one case.

Suppose that $g$ satisfies the inequalities (\ref{eq:jth_inequality}) for
each $1\leq j \leq r$.
Write $g=bk$, with $b\in B(\fbu)$ and $k\in G(\Oh)$.
Then $b$ can be written as $lu$, where $l\in L_{\{i\}}(\fbu)\cap B(\fbu)$
and
$u\in
U_{\{i\}}(\fbu)$.
Since $U_{\{i\}}$ fixes $\mathbf{v}_i$, the inequalities
(\ref{eq:jth_inequality}) imply that
\begin{equation}
\label{eq:new}
\| \mathbf{v}_i\rho_i(l)\|\leq \|\mathbf{v}_i\rho_i(\gamma l)\| \quad\mbox{ for
all }
\gamma \in L_{\{i\}}(R).
\end{equation}
From the rank one case, $l=\gamma \pi^\eta k$ for some $\gamma \in
L_{\{i\}}(R)$, $k\in L_{\{i\}}(\Oh)$ and $\eta \in X_*(T)$ such that
$|\alpha_i(\pi^\eta)|\geq 1$. $\rho_i(\gamma)$ maps $\mathbf{v}_i$ into
$V(R)$. From Lemma \ref{lemma:fundamental_inequality} it follows
that
\begin{equation*}
\|\mathbf{v}_i \rho_i(l)\| \geq \|\mathbf{v}_i \rho_i(\pi^\eta)\|.
\end{equation*}
Equation~(\ref{eq:new}) implies that the above must be an equality.
This forces $\gamma \in L_{\{i\}}(R)\cap P_i(R)$, and hence also $k \in
L_{\{i\}}(\Oh)\cap P_i(\Oh)$.
Write $b=tn$ with $t\in T(\fbu)$ and $n\in N(\fbu)$. Then viewing
$\alpha_i$ as a rational character of $B(\fbu)$ that is trivial on
$N(\fbu)$, we have
\begin{equation*}
|\alpha_i(t)|=|\alpha_i(l)|=|\alpha_i(\pi^\eta)|\geq 1.
\end{equation*}
Repeating this argument for each $i$ completes the proof of Lemma
\ref{lemma:otherthing}.\hfill \qed
\end{proof}

\section{Local reduction theory}\label{sec:local_proof}

In order to prove the existence part of Theorem \ref{theorem:local}, it
suffices to show that every element $g$ in $\Omega_G$ may be written as
$g=\gamma \pi^{\eta}k$, where $\gamma\in \Gamma$, $\eta\in X_*(T)$ is
antidominant and $k\in G(\Oh)$. To this end, we may assume (using the
Iwasawa decomposition) that we are given $g\in \Omega_G$, with $g=tn$,
with $t\in T(\fbu)$ and $n\in N(\fbu)$. Since $g$, and hence $t$, is in
$\Omega_G$, $|\alpha_i(t)|\geq 1$, so that $\alpha_i(t)^{-1}\in \Oh$,
for $i=1,\ldots,r$. For each root $\alpha\in \Phi(G,T)$, let $U_\alpha$
denote the corresponding root subgroup. Fix an isomorphism
$u_\alpha:\mathbf{G}_a \to U_\alpha$ defined over $\Fq$. Then for $x\in
\fbu$, we have
\begin{equation*}
tu_\alpha(x) = (tu_\alpha(x)t^{-1})t= u_\alpha(\alpha(t)x)t.
\end{equation*}
Therefore, if we write $\alpha(t)x=P+h$, where $P\in R$ and $h\in \Oh$,
then
\begin{equation*}
tu_\alpha(x) = t u_\alpha(\alpha(t)^{-1}P)u_\alpha(\alpha(t)^{-1}h) =
u_\alpha(P)t u_\alpha(\alpha(t)^{-1}h).
\end{equation*}
Given two positive roots $\alpha$ and $\beta$, the commutator
$[U_{\alpha},U_{\beta}]$ is contained in the product of root subgroups
$U_{\alpha'}$ where the $\alpha'$ are roots which can be written as
positive linear combinations of $\alpha$ and $\beta$ and are distinct
from either $\alpha$ or $\beta$. Moreover, we may enumerate the positive
roots as $\beta_1,\beta_2,\ldots$ so that if $j>i$, then $\beta_i$ cannot be written as a sum of $\beta_j$ and any other positive roots.

Write $n$ as $\prod_{i}u_{\beta_i}(x_i)$. Then
\begin{equation*}
tn = tu_{\beta_1}(x_1) \prod_{i>1} u_{\beta_i}(x_i).
\end{equation*}
If we write $\beta_1(t)x_1=P_1+h_1$, where $P_1\in \Fq[\pi^{-1}]$ and
$h\in \Oh$, then
\begin{equation*}
tn = u_{\beta_1}(P_1)tu_{\beta_1}(\beta_1(t)^{-1}h_1)\prod_{i>1}
u_{\beta_i}(x_i).
\end{equation*}
Since $u_{\beta_1}(P_1)\in \Gamma$, $\beta_1(t)^{-1}\in \Oh$, and the
image of $u_{\beta_1}$ normalizes all the subsequent root subgroups
whose elements appear in the above expression, we may assume for the
purpose of proving Theorem \ref{theorem:local}, that
\begin{equation*}
tn = t\prod_{i>1} u_{\beta_i}(x'_i),
\end{equation*}
for $x'_i\in \fbu$. We may continue in this manner to reduce $tn$ to
$t$. It is then easy to see (using the decomposition $\fbu^\times =
\pi^{\mathbf{Z}} \Oh^\times$) that $t$ may be replaced by $\pi^\eta$ for
$\eta \in X_*(T)$. Since $|\alpha_i(\pi^\eta)|\geq 1$, it follows that
$\eta$ is antidominant, proving the existence part of Theorem
\ref{theorem:local}.

We now prove the uniqueness part of Theorem \ref{theorem:local}. In
order to do this, it suffices to show that if $\eta$ and $\nu$ are two
dominant co-weights, and $\pi^\nu = \gamma \pi^\eta k$ for some $\gamma
\in \Gamma$ and $k\in G(\Oh)$, then $\nu=\eta$. Since the weights
$\mu_1,\ldots, \mu_r$ corresponding to the fundamental representations
in Theorem \ref{theorem:godement} generate the vector space
$X^*(T)\otimes \mathbf{Q}$, it suffices to show that $\langle \mu_i, \nu
\rangle = \langle \mu_i, \eta \rangle$ for each $i$. In order to do
this, we need the following:
\setcounter{lem}{23}
\begin{lem}
\label{lemma:fundamental_inequality}
For any non-zero vector $\mathbf{v}\in V_i(\fbu)$ and any antidominant
co-weight $\mu\in X_*(T),$
\begin{equation*}
\frac{\|\mathbf{v}\rho_i(\pi^\mu)\|}{\|\mathbf{v}\|}\geq
\frac{\|\mathbf{v}_i\rho_i(\pi^\mu)\|}{\|\mathbf{v}_i\|}.
\end{equation*}
\end{lem}
\begin{proof}
Since $T$ is an $\Fq$-split torus and $\rho_i$ is defined over $\Fq$,
$V$ has a decomposition (over $\Fq$) into root subspaces
\begin{equation*}
V = \bigoplus_\lambda V_\lambda,
\end{equation*}
where $T$ acts on $V_\lambda$ by the character $\lambda:T\to \mathbf{G}_m$.
It is easy to see that $\mu_i$ is the lowest weight of $T$ occurring in $(\rho_i,V_i)$, so that $\langle \mu_i,\mu \rangle\geq \langle \lambda,\mu \rangle$ for any weight $\lambda$ of $T$ occurring in $(\rho_i,V_i)$ and any antidominant co-weight $\mu$.
Given any vector $\mathbf{v}\in V(F_\bu)$, we may write
\begin{equation*}
\mathbf{v}=\sum x_j\mathbf{u}_j,
\end{equation*}
where $x_j\in F_\bu$ and $\mathbf{u}_j\in V_{\lambda_j}(\Fq)$ for each
$j$ and the $\lambda_j$s are not necessarily distinct. Thus
\begin{align*}
\| \mathbf{v}\rho_i(\pi^\mu) \| &=\left\| \sum \lambda_j(\pi^\mu)x_j \mathbf{u}_j \right\|\\[.3pc]
& = \sup_j\{ | \lambda_j(\pi^\mu)x_j| \}\\[.3pc]
& = \sup_j\{ q^{-\langle \lambda_j, \mu \rangle}|x_j|\}\\[.3pc]
& \geq q^{-\langle \mu_i,\mu \rangle} \sup_j\{|x_j|\}\\[.3pc]
& = \| \mathbf{v}_i \rho_i(\pi^\mu)\|\; \|\mathbf{v}\|.
\end{align*}
Since $\|\mathbf{v}_i\|=1$, this completes the proof of Lemma
\ref{lemma:fundamental_inequality}.\hfill \qed
\end{proof}

Lemma \ref{lemma:fundamental_inequality} allows us to compare $\langle
\mu_i,\nu\rangle$ and $\langle \mu_i , \eta \rangle$:
\begin{align*}
q^{-\langle \mu_i,\eta\rangle}& = \frac{\|\mathbf{v}_i \rho_i(\pi^\eta)\|}{\|\mathbf{v}_i\|}\\[.3pc]
& \leq \frac{\|\mathbf{v}_i\rho_i(\gamma \pi^\eta)\|}{\|\mathbf{v}_i\rho_i(\gamma)\|}\\[.3pc]
& \leq \frac{\| \mathbf{v}_i\rho_i(\gamma \pi^\eta)\|}{\|\mathbf{v}_i\|}\\[.3pc]
& = \frac{\|\mathbf{v}_i\rho_1(\pi^\nu)\|}{\|\mathbf{v}_i\|}\\[.3pc]
&= q^{-\langle \mu_i, \nu \rangle}.
\end{align*}
The first inequality is Lemma \ref{lemma:fundamental_inequality} applied
to $\mathbf{v}=\mathbf{v}_i\rho_i(\gamma)$. The second inequality
follows from Lemma \ref{lemma:boundone} with
$\mathbf{x}=\mathbf{v}_i\rho_i(\gamma)$. Interchanging the roles of
$\eta$ and $\nu$ in the above arguments shows that $\langle \mu_i,\eta
\rangle=\langle \mu_i, \nu \rangle$ for each $i$. This completes the
proof of the uniqueness part of the assertion of Theorem
\ref{theorem:local}.

\section{Global reduction theory}
\label{sec:global_proof}

If $g= (g_v)_v$ is an element of $G(\A)$ then, since $g_v \in G(\Oh_v)$
for all but finitely many places $v$ of $F$, we may assume, for the
purpose of proving Theorem \ref{theorem:global_reduction} that $g$ is a
finite product $g=g_\infty g_{v_1} g_{v_2} \ldots g_{v_k}$, with
$g_\infty \in G(F_\infty)$ and $g_{v_j}\in G(F_{v_j})$, $v_j\neq\infty$,
for $1\leq j \leq k$. By Theorem \ref{theorem:local}, there is a
decomposition
\begin{equation*}
g_{v_k} = \gamma_k \pi_{v_k}^{\eta_k} \kappa_k,
\end{equation*}
where $\gamma_k \in G(\Fq[\pi_{v_k}^{-1}])$, $\eta_k \in X_*(T)$, and
$\kappa_k \in G(\Oh_{v_k})$. Now $\gamma_k$ and $\pi_{v_k}^{\eta_k}$ are
both contained in $G(F)$ and in $G(\Oh_v)$ for all $v\neq \infty$.
Therefore, by multiplying $g$ on the left by
$\pi_{v_k}^{-\eta_k}\gamma^{-1}$ we get an element of the subset
\begin{equation*}
G(F_\infty)\times \prod_{j=1}^{k-1}G(F_{v_j}) \times \prod_{{\rm all\ other}\ v} G(\Oh_v)
\end{equation*}
of $G(\A)$.

We have now reduced $g$ to an element with non-trivial entries  only
at most $k-1$ places and $\infty$. We may continue in this manner until
the entries at all places except $\infty$ are trivial. Finally, the use
of Theorem \ref{theorem:local} to $v=\infty$ gives us a representative
each double coset of type asserted by Theorem
\ref{theorem:global_reduction}.

The uniqueness part of the theorem follows from the corresponding
assertion in the local situation, because two elements $g$ and $h$ of
$G(F_\infty)$ lie in the same double coset if and only if $g=\gamma h
k$, with $\gamma \in G(\Fq[t])$ and $k\in G(\Oh_\infty)$.

\section*{Acknowledgements}

We thank Robert Kottwitz for his guidance and R Narasimhan and
Dipendra Prasad for some historical and bibliographical information. We
thank the University of Chicago and the Centre de Recherches
Math\'ematiques for their support during the preparation of this article.

\end{document}